# Fermat's Last Theorem on Topological Fields

Vinod Kumar.P.B[1] & K.Babu Joseph[2]

Department of Mathematics

Rajagiri School of Engineering and Technology

Rajagiri Valley.P.O, Kochi – 682039

Kerala, India.

ABSTRACT:

Fermat's last theorem – even though it is a number theoretic result we prove that the result depends on the topological as well as the field structure of the underlying space. Geometry of the space can be interpreted as a topological property. A change in geometry or rather in topology may violate FLT. Andrew Wiles has proved that FLT is true. But we shall prove that it is because of the Euclidean structure of $R^n$. So "FLT is true" is not a globally valid statement, rather it is relative that it depends on the topological field structure of the space.

2000 AMS classification: 54E99, 11T99.

1 vinod_kumar@rajagiritech.ac.in

2 babu_65@hotmail.com

1. **Introduction**

In 1637 Pierre de Fermat conjectured that there are no integer solutions for $x^n + y^n = z^n$, if $n \succ 2$. If n=2 it can be interpreted as Pythagorus theorem. Fermat, motivated by the work of Diophantus argued that he know the proof but did not publish it. The case for n= 4 was proved by Fermat himself. For n=3, proof was given by Leonhard Euler in 1770. For n=5, Dirichlet and Legendre proved it in 1825. For n=7, Gabriel Lamme proved it in 1890. If the numbers x,y,z are a solution to FLT for n = 6, then the numbers $x^2$, $y^2$, $z^2$ are a solution for n = 3, contradicting FLT for n = 3. So generally if we had proven FLT to any power k, then the theorem is valid for all the multiples of k. So it is enough to prove FLT just for powers that are primes. The proof of Andrew Wiles, valid for $n \succ 2$ is a consequence of two results.

(i)  if there is a solution to FLT, then the elliptic curve defined by the equation $Y^2 = (X - x^n)(Y - y^n)$ is semistable, but not modular.

(ii)  All semistable elliptic curves with rational coefficients are modular.

In 1986, Ken Ribet and Jean Piere Serie proved result (i) and A.Wiles with the help of R.taylor proved result(ii) in 1995.In fact, result(ii) was conjectured by Y.taniyama in 1955 and in more general form by G.Shimura and A. Weil. So it is clear that Wiles proof depends on the concept of elliptic curves.Elliptic curves have the geometrical structure. The geometry will vary as topology. So we claim that it is not true that "FLT is true" in general. FLT is true in $R^n$ with the usual topology. In this paper a relation between topology, field theory , number theory and FLT is established.

Section 1 deals with some results in topological fields and use some of them in establishing our claim. In section 3, we discuss FLT on finite fields. We will prove that FLT is true in finite field $F_p$ only if n = (p-1) or (p-1)/2. In section 4, the geometry of topological fields giving special importance to elliptic curves is discussed. Our claim will be proved in section5, in which we discuss FLT on topological fields.

1. **Topological Fields**

Definitions

(i) Let (F,*,.) be a Field and T be a topology on F. Then (F,*,.,T) is called a topological Field.

There are topological fields which look like as same if we view from different angles, for example- $R^n$.and there are topological fields which is not same at all points for example – the figure **8.** View of '8' from its saddle point is different from its view from any other point. In the first case we say that the space is homogenous.

(ii) A topological space F is homogenous if given any x,y $\in$ F, there exist f: F $\to$ F, a homeomorphism on to F such that f(x) = y.

(iii)    In (F,*,.,T) neighborhood of 0 is called additive nucleius.

We prove the following results in this section.

(i)    Every topological field is homogenous.
(ii)   Every topological field is Hausdorff.

Proofs:

(i)    : Let x, y $\in$ F. Let $L_a(b)$ = a.b.

Then $L_{yx}^{-1}(x)$ = y. So given x,y $\in$ F, there exists $L_{yx}^{-1}$: F $\to$ F, which is a homeomorphism.

(ii)   Let x, y $\in$ F.

$\Rightarrow x^{-1}.y \neq 0$. $\Rightarrow$ X-$\{x^{-1}.y\}$ is a neighbourhood of 0. $\Rightarrow$ there is a symmetric neighbourhood V $\subseteq$ X-$\{x^{-1}.y\}$. $\Rightarrow$ x.V and y.V are nbds. of x and y. $\Rightarrow$ x.V $\cap$ y.V $\neq \phi$.

So there exist $v_1, v_2 \in$ V such that $v_1 v_2^{-1} \in$ V. Which is a contradiction. Hence F is Hausdorff.

2. **Finite Fields**

Denote $\mathbb{F}[p^n]$, where $p$ is a prime $n$ is a positive integer $\mathbb{F}[p^n]$ has $p^n$ elements

For n=1, fields are of the form

F(p) = {0,1,2,3,...,p-1}
multiplication and addition are usual operations, except multiples of p should be left out of the set (modulo p).
$$F(5) = \{0,1,2,3,4)$$
All of the multiplication in the example is mod 5 because $p=5$ .

Addition on $\mathbb{F}_{[4]}$ is addition of polynomials. $(a+bx)+(c+dx) = (a+c)+x(b+d)$ What is the additive identity?

$(a+bx)+(c+dx) = (a+bx) = (a+c)+x(b+d)$

a+c= 0 c= 0

b+d= b d= 0

∴ 0 is still the additive identity.

(a+bx) has the additive inverse -a-bx.

usual multiplication of polynomials *except* need to use $x^2+x+1=0$.

$(a+bx)*(c+dx) = ac = x(bc+ad)+x^2(bd)$

this seems to be a problem because $\mathbb{F}_{[4]}$ only has linear polynomials but we ended up with a quadratic one.

When constructing the elements of $\mathbb{F}_{[2]}$ remember had multiplication modulo p = 2, so $\mathbb{F}_{[4]}$ has polynomial multiplication modulo $x^2+x+1$. Also $\mathbb{F}_{[2]}$ sits inside of $\mathbb{F}_{[4]}$ as the "constant" polynomials.

⇒ Remember when making these tables that each element will only show up once in any column or row. This is because we want to show that each element has a multiplicative inverse, and that there are no zero divisors.

As a vector space, $\mathbb{F}_{[9]} = \mathbb{F}_{[3^2]} = \{a+bx | a,b \in \mathbb{F}_{[3]}\}$

To find the multiplication table we need a monic-quadratic that has no zeros in $\mathbb{F}_{[3]}$

A monic-quadratic will have a coefficient of 1 on the highest degree term.

<u>Try</u> $1x^2+0x+1$, f(0)= 1, f(1)=2 f(2)=2 (no zeros!)

Note: The constant term must to be non zero, because otherwise 0 is a zero.

Also note that we only need one polynomial that works.

From the above explanation, it is clear that in $F_p$, FLT is true only if n – (p-1) or (p-1)/2.

If we change the metric on R2, its geometry is different and hence FLT is no more true.

For example if we consider the taxi cab metric "d" on $R^2$,

$d((x_1,y_1), (x_2,y_2)) = |x_1-x_2| + |y_1-y_2|$,

$[d((0,0),(1,0))]^2 + [d((0,0),(0,1))]^2 \neq [d((1,0),(0,1))]^2$.

Geometric shapes do not remain as it is if we change the topology of $R^n$.

We will discuss the geometry of elliptic curves in the next section

### 3. Elliptic curves

Elliptic curve is simply the locus of points in the x-y plane that satisfy an algebraic equation of the form $Y^2 = (X - x^n)(Y - y^n)$ (with some additional minor technical conditions). This is deliberately vague as to what sort of values x and y represent. In the most elementary case, they are real numbers, in which case the elliptic curve is easily graphed in the usual Cartesian plane. But the theory is much richer when x and y may be any complex numbers (in **C**). And for arithmetic purposes, x and y may lie in some other field, such as the rational numbers **Q** or a finite field **F**.

So an elliptic curve is an object that is easily definable with simple high school algebra. Its amazing fruitfulness as an object of investigation may well depend on this simplicity, which makes possible the study of a number of much more sophisticated mathematical objects that can be defined in terms of elliptic curves.

It is very natural to work with curves in the complex numbers, since **C** is the *algebraic closure* of the real numbers. That is, it is the smallest algebraically closed field that contains the roots of all possible polynominals with coefficients in **R**. Being algebraically closed means that **C** contains the roots of all polynomials with coefficients in **C** itself. It's natural to work with a curve in an algebraically closed field, since then the curve is as "full" as possible.

The case of elliptic curves in the complex numbers is especially interesting, not only because of the algebraic completeness of **C**, but also because of the rich analytic theory that exists for complex functions. In particular, the equation of an elliptic curve defines y as an "algebraic function" of x. For every algebraic function, it is possible to construct a specific surface such that the function is "single-valued" on the surface as a domain of definition. It turns out that an elliptic curve, defined as a locus of points, is also the Riemann surface associated with the algebraic function defined by the equation.

So an elliptic curve is a Riemann surface. In fact, it is of a special type: a compact Riemann surface of genus 1. And not only that, but the converse is also true: every compact Riemann surface of genus 1 is an elliptic curve. In other words, elliptic curves over the complex numbers represent exactly the "simplest" sorts of compact Riemann surfaces with non-zero genus. Topologically, the genus counts the number of "holes" in a surface. A surface with one hole is a torus.

This topological equivalence of an elliptic curve with a torus is actually given by an explicit mapping involving the Weierstrass $\wp$-function and its first derivative. This mapping is, in effect, a parameterization of the elliptic curve by points in a "fundamental parallelogram" in the complex plane.

The topological space that results from identifying opposite sides of a period parallelogram is called a complex torus. The fundamental periods that define the parallellogram generate a *lattice* in **C** consisting of all sums of integral multiples of $\omega_1$ and $\omega_2$. If L denotes the lattice, then L = $\mathbf{Z}\omega_1 \oplus \mathbf{Z}\omega_2$. The complex torus can then be described as **C**/L. What this all means, therefore, is that a (complex) torus is the "natural" domain of definition of the $\wp$-function, or any doubly periodic complex function.

. The Fermat equation is the prime example. In general, an elliptic curve has the form $y^2 = Ax^3 + Bx^2 + Cx + D$, but for considering arithmetical questions, it is natural to restrict our attention to the case where A, B, C, D are all rational. This assumption will usually be in effect when we are considering properties of elliptic curves involving arithmetical questions (as opposed to their more general analytic properties). If all coefficients are rational, the elliptic curve is said to be *defined over* **Q**. The all-important Taniyama-Shimura conjecture concerns only elliptic curves defined over **Q**.

The fact that any elliptic curve (not necessarily defined over **Q**) has an abelian group structure means that we can learn a lot about it by studying various of its subgroups. For considering arithmetical (i. e. number theoretic) questions, we restrict our attention to curves defined over **Q**. In that case, there are several interesting subgroups we can consider.

The first is the group of all points on the curve E which have an order that divides m for some particular integer m. That is, m "times" such a point is the identity element. Such points are called "m-division points", and the subgroup they make up is denoted E[m]. The reason for the name is that any point in E[m] generates a cyclic subroup of E (and E[m]) whose order divides m. If the order is actually m, then the points in the cyclic group generated by the point divide E into m segments.

It isn't necessarily the case that the coordinates of a point in E[m] have integral or rational coordinates. However, the coordinates will be algebraic numbers (i. e., roots of an algebraic equation with coefficients in **Q**). It's relatively easy to show that as an abstract group E[m] is just the direct sum of two cyclic groups of order m, i. e. $\mathbf{Z}/m\mathbf{Z} \oplus \mathbf{Z}/m\mathbf{Z}$, so its order is $m^2$. We shall see later that its real interest lies in the fact that we can construct representations of other groups of transformations that act on E[m]. Such representations will consist of 2x2 matrices with integral entries, i. e. elements of $GL_2(\mathbf{Z})$.

Another interesting subgroup of E is the set of all points whose coordinates are rational. Such points are said to be *rational points*. If the curve is defined over **Q**, then it is a simple fact that the set of all rational points (if there are any) is a subgroup.

The definition of L(E,s) will be made based on details about a series of other groups connected with E. These arise by considering E as an elliptic curve over the finite fields $\mathbf{F}_p$. This is the same as taking the original equation and reducing the coefficients mod p. If the equation of E has rational but non-integral coefficients, we would need to assume none of their denominators are divisible by p, so we might as well assume all coefficients to be integral to begin with (since if the denominators are prime to p they have inverses mod p). Further, the definition of an elliptic curve requires that there are no repeated roots of the polynomial in x, and this may fail to be true when reducing mod p for some primes. Such primes are said to have "bad reduction". There will be only a finite number of these for any particular curve (they will divide the discriminant), but they have to be dealt with specially.

For any prime p where E has good reduction, we can consider the elliptic curve $E(\mathbf{F}_p)$ over $\mathbf{F}_p$. Since $\mathbf{F}_p$ is finite, there are only a finite number of points on $E(\mathbf{F}_p)$, so it is a finite group. The order of this group, $\#(E(\mathbf{F}_p))$, turns out to be a very important number.

There is a general approach in number theory of trying to deal with "global" problems, such investigating the structure of E(**Q**), by looking at a closely related "local" problem mod p for all primes p. This is why we are interested in $E(\mathbf{F}_p)$. In particular, if $E(\mathbf{F}_p)$ is "large" for most p, we would expect E(**Q**) to be large too.

We will see that the numbers $\#(E(\mathbf{F}_p))$ are studied by relating them to coefficients of the Dirichlet series of L(E,s), the L-function of E.

The most important fact about the minimal discriminant is that the primes which divide it are precisely the ones at which the curve has bad reduction. In other words, except for those primes, the reduced curve is an elliptic curve over $\mathbf{F}_p$.

There is still another invariant of an elliptic curve E, called its *conductor*, and often denoted simply by N. The exact definition is rather technical, but basically the conductor is, like the minimal discriminant, a product of primes at which the curve has bad reduction. Recall that E has bad reduction when it has a singularity modulo p. The type of singularity determines the power of p that occurs in the conductor. If the singularity is a "node", corresponding to a double root of the polynomial, the curve is said to have "multiplicative reduction" and p occurs to the first power in the conductor. If the singularity is a "cusp", corresponding to a triple root, E is said to have "additive reduction", and p occurs in the conductor with a power of 2 or more.

If the conductor of E is N, then it will turn out that N is the "level" of certain functions called modular forms (not yet defined) with which, according to the Taniyama-Shimura conjecture, E is intimately connected.

If N is square-free, then all cases of bad reduction are of the multiplicative type. An elliptic curve of this sort is called *semistable*. It is for elliptic curves of this sort that Wiles proved the Taniyama-Shimura conjecture.

**Theorem 1** - *Any finite field with characteristic p has $p^n$ elements for some positive integer n.*

Proof: Let L be the finite field and K the prime subfield of L. The vector space of L over K is of some finite dimension, say n, and there exists a basis $\þ_1, \þ_2, \ldots, \þ_n$ of L over K. Since every element of L can be expressed uniquely as a linear combination of the $\þ_i$ over K, i.e., every a in L can be written as, a = Sum $\ß_i \þ_i$, with $\ß_i$ in K, and since K has p elements, L must have $p^n$ elements. ¶

This theorem, while it does restrict the size of a finite field, does not say that one will exist for a particular power of a prime, nor does it specify how many finite fields can exist of a particular order. The answers to these questions can be deduced from the following theorem.

**Theorem 2** - *Let L be a field with characteristic p and prime subfield K. Then L is the splitting field for*

$$f(x) = x^{p^n} - x$$

*iff L has $p^n$ elements.*

Proof: Suppose that L is the splitting field for $f(x) = x^{p^n} - x$ over K. Since $(f(x), f'(x)) = 1$, the roots of $f(x)$ are distinct and so L has at least $p^n$ elements. Consider the subset

$E = \{ \text{þ in L} \mid \text{þ}^{p^n} = \text{þ} \}$

of L. Clearly E contains $p^n$ elements since it consists of the roots of $f(x)$. Suppose that þ,ß in E; then $(\text{þß})^{p^n} = (\text{þ})^{p^n} (\text{ß})^{p^n} = \text{þß}$ and hence, þß in E. Also,

$$(\text{þ}+\beta)^{p^n} = \sum \binom{p^n}{i} \text{þ}^i \beta^{p^n-i} = \text{þ}^{p^n} + \beta^{p^n} = \text{þ} + \beta$$

since $p \mid C(p^n, i)$ for $0 < i < p^n$, and hence, $(\text{þ} + \text{ß})$ in E. The existence of additive and multiplicative inverses is easy to show, so E is a subfield of L and also a splitting field for $f(x)$. Thus by Thm II.1.1 $E = L$ and L contains $p^n$ elements.

Suppose now that L contains $p^n$ elements. The multiplicative group of L, which we will denote by L*, forms a group of order $p^n - 1$ and hence the order of any element of L* divides $p^n - 1$. Thus $\text{þ}^{p^n} = \text{þ}$ for all þ in L* and the relation is trivially true for þ = 0. Thus $f(x)$ splits in L. ¶

Now for some important corollaries.

**Corollary.5** - *For any prime p and integer n,* **GF**$(p^n)$ *exists*.

Proof: By Thm II.1.1 the splitting field exists and by Cor II.2.4 it is GF$(p^n)$. ¶

The following important theorem is useful in establishing the subfield structure of the Galois Fields among other things.

**Theorem 6** - **GF**$(p^n)$* *is cyclic*.

Proof: The multiplicative group **GF**$(p^n)$* is, by definition, abelian and of order $p^n - 1$. If $p^n - 1 = p_1^{e_1} \ldots p_k^{e_k}$, then, factoring **GF**$(p^n)$* into a direct product of its Sylow subgroups, we have **GF**$(p^n)$* = $S(p_1) \times \ldots \times S(p_k)$ where $S(p_i)$ is the Sylow subgroup of order $(p_i)^{e_i}$. The order of every element in $S(p_i)$ is a power of $p_i$ and let $a_i$ in $S(p_i)$ have the maximal order, say $(p_i)^{e'_i}$, $e'_i <= e_i$, for $i = 1,...,k$. Since $(p_i, p_j) = 1$, i not equal j, the element $a = a_1 a_2 \ldots a_k$ has maximal order $m = (p_1)^{e'_1} \ldots (p_k)^{e'_k}$ in $GF(p^n)*$. Furthermore every element of $GF(p^n)*$ satisfies the polynomial $x^m -1$, implying that $m >= p^n -1$. Since a in $GF(p^n)*$ has order m, m divides $p^n -1$ and so, $m = p^n -1$. Thus the element a is a generator and $GF(p^n)*$ is cyclic. ¶

A generator of $GF(p^n)*$ is called a *primitive element* of $GF(p^n)$.

The following theorem has some useful consequences.

**Theorem 7** - *Over any field K, $(x^m - 1) | (x^n - 1)$ iff $m | n$.*

Proof: If $n = qm + r$, with $r < m$, then by direct computation

$$x^n - 1 = x^r \left( \sum_{i=0}^{q-1} x^{im} \right) (x^m - 1) + (x^r - 1).$$

It follows that $(x^m - 1) | (x^n - 1)$ iff $x^r - 1 = 0$, i.e., $r = 0$. ¶

**Corollary 8** - *For any prime integer p, $(p^m - 1) | (p^n - 1)$ iff $m | n$.*

Proof: Basically the same as that of the theorem, do for homework.

**Theorem 9** - $GF(p^m)$ *is a subfield of* $GF(p^n)$ *iff $m | n$.*

Proof: Suppose $GF(p^m)$ is a subfield of $GF(p^n)$; then $GF(p^n)$ may be interpreted as a vector space over $GF(p^m)$ with dimension, say, k. Hence, $p^n = p^{km}$ and $m | n$.

Now suppose $m | n$, which from the previous theorem and its corollary implies that $(x^{p^m - 1} - 1) | (x^{p^n - 1} - 1)$. Thus every zero of $x^{p^m} - x$ that is in $GF(p^m)$ is also a zero of $x^{p^n} - x$ and hence in $GF(p^n)$. It follows that $GF(p^m)$ is contained in $GF(p^n)$. Notice that there is precisely one subfield of $GF(p^n)$ of order $p^m$, otherwise $x^{p^m} - x$ would have more than $p^m$ roots. ¶

Although we will not prove it, the automorphism group of a finite field is cyclic. The standard generator of this group is the so-called *Frobenius automorphism* defined for a finite field of characteristic p as the map $x \longrightarrow x^p$ for all x in $GF(p^n)$.

4. **Fermat's Last Theorem**

From the discussions and theorems in Section2, 3, and 4 it is clear that FLT is based on the theory of elliptic curves which in turn depends on the topology of the space. So precisely from Theorem 2 and 3 in section 3 ; and from all the 9 theorems in section 4 our claim is established.

**Acknowledgement:**The authors thank the management of the institution for the support and inspiration.

************